\documentclass [11pt] {article}
\usepackage{amssymb}
\usepackage{latexsym}
\usepackage{amsmath}
\title{Some Remarks on\\
       \mbox{$\mathfrak g$}-invariant Fedosov Star Products and\\
       Quantum Momentum Mappings}
\author{{\bf Michael Frank M\"uller
\thanks{Michael.Mueller@math.uni-mannheim.de}}\\[3mm]
Fakult\"at f\"ur Mathematik und Informatik\\Universit\"at
Mannheim\\A5 Geb\"aude C\\D-68131 Mannheim\\ Germany\\[3mm] {\bf
Nikolai Neumaier
\thanks{Nikolai.Neumaier@physik.uni-freiburg.de}}\\[3mm]
Fakult\"at f\"ur Mathematik und Physik\\Universit\"at
Freiburg\\Hermann-Herder-Stra\ss e 3\\D-79104 Freiburg i.~Br.\\
Germany\\[3mm]}
\date{January 2003\\[3mm]FR-THEP-2003/01
\\[3mm]Mannheimer Manuskripte 268}

\newcommand {\Cinf} [1] {\mathcal C^\infty ({#1})}
\newcommand {\Ginf} [1] {\Gamma^\infty ({#1})}
\newcommand {\JN} [1] {J_0({#1})}
\newcommand {\JP} [1] {J_+({#1})}
\newcommand {\Jbold} [1] {J({#1})}
\newcommand {\symp} {\Gamma^\infty_{\mbox{\rm\tiny symp}}(TM)}

\newcommand {\id} {{\mathrm{id}}}
\newcommand {\Lie} {{\mathcal L}}
\newcommand {\ad} {{\mathrm{ad}}}
\newcommand {\W} {\mathcal W}
\newcommand {\WL} {\mbox{$\mathcal W \! \otimes \! \Lambda$}}
\newcommand {\degs} {{\rm deg}_s}
\newcommand {\dega} {{\rm deg}_a}
\newcommand {\degnu} {{\rm deg}_\nu}
\newcommand {\Deg} {{\rm Deg}}
\renewcommand {\d} {\mathrm{d}}

\newenvironment {PROOF}{\small {\sc Proof:}}{{\hspace*{\fill}
                       $\square$}}

\newtheorem {LEMMA} {Lemma} [section]

\newtheorem {PROPOSITION} [LEMMA] {Proposition}
\newtheorem {THEOREM} [LEMMA] {Theorem}
\newtheorem {COROLLARY} [LEMMA] {Corollary}
\newtheorem {DEFINITION}[LEMMA] {Definition}
\newtheorem {REMARK}[LEMMA] {Remark}
\newtheorem {DEDUCTION}[LEMMA] {Deduction}
\newtheorem {EXAMPLE} [LEMMA] {Example}
\begin{document}
\maketitle
\begin{abstract}
In these notes we consider the usual Fedosov star product on a
symplectic manifold $(M,\omega)$ emanating from the fibrewise Weyl
product $\circ$, a symplectic torsion free connection $\nabla$ on
$M$, a formal series $\Omega \in \nu Z^2_{\mbox{\rm\tiny
dR}}(M)[[\nu]]$ of closed two-forms on $M$ and a certain formal
series $s$ of symmetric contravariant tensor fields on $M$. For a
given symplectic vector field $X$ on $M$ we derive necessary and
sufficient conditions for the triple $(\nabla,\Omega,s)$
determining the star product $*$ on which the Lie derivative
$\Lie_X$ with respect to $X$ is a derivation of $*$. Moreover, we
also give additional conditions on which $\Lie_X$ is even a
quasi-inner derivation. Using these results we find necessary and
sufficient criteria for a Fedosov star product to be $\mathfrak
g$-invariant and to admit a quantum Hamiltonian. Finally, supposing
the existence of a quantum Hamiltonian, we present a cohomological
condition on $\Omega$ that is equivalent to the existence of a
quantum momentum mapping. In particular, our results show that the
existence of a classical momentum mapping in general does not imply
the existence of a quantum momentum mapping.
\end{abstract}
\clearpage
\tableofcontents
\section{Introduction}
\label{IntroSec}
The concept of deformation quantization as introduced in the
pioneering articles \cite{BayFla78} by Bayen, Flato, Fr\o nsdal,
Lichnerowicz and Sternheimer has proved to be an extremely useful
framework for the problem of quantization: the question of
existence of star products $\star$ (i.e. formal, associative
deformations of the classical Poisson algebra of complex-valued
functions $\Cinf{M}$ on a symplectic or more generally, on a
Poisson manifold $M$, such that in the first order of the formal
parameter $\nu$ the commutator of the star product yields the
Poisson bracket) has been answered positively by DeWilde and
Lecomte \cite{DeWLec83b}, Fedosov \cite{Fed94}, Omori, Maeda and
Yoshioka \cite{OmoMaeYos91} in the case of a symplectic phase space
as well as by Kontsevich \cite{Kon97} in the more general case of a
Poisson manifold. Moreover, star products have been classified up
to equivalence in terms of geometrical data of the phase space by
Nest and Tsygan \cite{NT95a}, Bertelson, Cahen and Gutt
\cite{BerCahGut97}, Weinstein and Xu \cite{WX97} on symplectic
manifolds and the classification on Poisson manifolds is due to
Kontsevich \cite{Kon97}. Comparisons between the different results
on classification and reviews can be found in articles of Deligne
\cite{Del95}, Gutt and Rawnsley \cite{GR99,Gut00}, Neumaier
\cite{Neu99} and Dito and Sternheimer \cite{DitSte02,Ste98}.

Already at the very beginning of the investigations of deformation
quantization various notions of invariance of star products with
respect to Lie group resp. Lie algebra actions were introduced and
discussed by Arnal, Cortet, Molin and Pinczon in \cite{ArnCor83}.
Later on it was Xu who systematically defined the notion of a
quantum momentum mapping for $\mathfrak g$-invariant star products
in the framework of deformation quantization in \cite{Xu98} that
naturally generalizes the concept of the momentum mapping in
Hamiltonian mechanics (cf. \cite{AbrMar85}) and computed the a
priori obstructions for its existence. Actually the notion of a
quantum momentum mapping has proved to be essential for the
formulation of the quantum mechanical analogue of the
Marsden-Weinstein reduction in deformation quantization as it was
studied by Fedosov in \cite{Fed98}, where it was shown that in some
sense `reduction commutes with quantization'. For the application
of the BRST quantization in deformation quantization as it was
introduced and discussed by Bordemann, Herbig and Waldmann in
\cite{BorHerWal00} the existence of a quantum momentum mapping also
turned out to be a major ingredient of the construction. For the
more special discussion of the example of reduction of star
products for $\mathbb C P^n$ as it was given by Bordemann,
Brischle, Emmrich and Waldmann in \cite{BorBriEmm96} and was
slightly generalized by Waldmann in \cite{Wal98} again the use of a
quantum momentum mapping the existence of which can be shown
explicitly in this case was the key ingredient of the
considerations.

Recently in \cite{Ham02} Hamachi has taken up afresh the question
under which preconditions the usual Fedosov star product admits a
quantum momentum mapping and he has given a condition in terms of
parts of the Fedosov derivation used to define the star product
which is assumed to be invariant with respect to a symplectic Lie
Group action on $M$.

In the present paper we want to generalize these results into two
directions: Firstly we drop the assumption of invariance of the
star product with respect to a Lie group action and replace it by
the somewhat weaker invariance with respect to the action of a Lie
algebra $\mathfrak g$. Secondly we make the conditions given in
\cite{Ham02} more precise and show that assuming that there is a
classical momentum mapping the question of existence of a quantum
momentum mapping relies on two cohomological conditions on the
formal series $\Omega \in \nu Z^2_{\mbox{\rm\tiny dR}}(M)[[\nu]]$
used to construct the $\mathfrak g$-invariant star product.

The paper is organized as follows: In Section \ref{PrelimSec} we
collect some notations and give a very short review of Fedosov's
construction. Here we also prove some technical details that enable
us to describe all derivations of the Fedosov star products in a
very convenient way which turns out to be very useful for the
further investigations. In Section \ref{SymVecDerSec} we consider
an arbitrary symplectic vector field on $M$ and give necessary and
sufficient conditions for the Lie derivative with respect to this
vector field to be a derivation of the star product $*$ under
consideration. Furthermore we can also specify additional
conditions guaranteeing that this derivation is even quasi-inner.
In Section \ref{ginvQMMapSec} we recall the definitions of
$\mathfrak g$-invariant star products, quantum Hamiltonians and
quantum momentum mappings from \cite{Xu98} and apply our result of
Section \ref{SymVecDerSec} to give criteria for the $\mathfrak
g$-invariance of a Fedosov star product. Finally, supposing that
the Lie algebra action is Hamiltonian and the Hamiltonian is
equivariant with respect to the coadjoint action of $\mathfrak g$
we moreover find conditions that permit a decision whether quantum
momentum mappings do exist. We conclude the paper with some remarks
on possible generalizations and further investigations.

{\bf Conventions:} By $\Cinf{M}$, we denote the complex-valued
smooth functions and similarly $\Ginf{T^*M}$ stands for the
complex-valued smooth one-forms et cetera. Moreover, we use
Einstein's summation convention in local expressions.

\section{Preliminaries}
\label{PrelimSec}
In this section we shall briefly recall the essentials of Fedosov's
construction of star products on a symplectic manifold
$(M,\omega)$. As we assume the reader to be familiar with this
construction we shall restrict to the very minimum to introduce our
notation (For more details we refer the reader to
\cite{Fed94,Fed96} and \cite[Sect. 2]{Neu99}, where we even used
the same notation). Defining
\begin{equation} \label{WeylAlgDef}
    \WL := \left({\mathsf X}_{s=0}^\infty
               \Gamma^\infty \left(\mbox{$\bigvee$}^s T^*M
               \otimes \mbox{$\bigwedge$}
               T^*M\right)\right)[[\nu]].
\end{equation}
it is obvious that $\WL$ becomes in a natural way an associative,
super-commutative algebra and the product is denoted by $\mu (a
\otimes b) = ab$ for $a,b\in \WL$ (By $\WL^k$ we denote the
elements of anti-symmetric degree $k$ and set $\W := \WL^0$.).
Besides this pointwise product the Poisson tensor $\Lambda$
corresponding to $\omega$ gives rise to another associative product
$\circ$ on $\WL$ by
\begin{equation} \label{FibProd}
    a \circ b = \mu \circ \exp\left(\frac{\nu}{2}\Lambda^{ij}
    i_s(\partial_i) \otimes i_s(\partial_j)\right)(a \otimes b),
\end{equation}
which is a deformation of $\mu$. Here $i_s(Y)$ denotes the
symmetric insertion of a vector field $Y\in \Ginf{TM}$ and
similarly $i_a(Y)$ shall be used to denote the anti-symmetric
insertion of a vector field. We set $\ad(a)b
:= [a,b]$ where the latter denotes the $\dega$-graded
super-commutator with respect to $\circ$. Denoting the obvious
degree-maps by $\degs$, $\dega$ and $\degnu = \nu \partial_\nu$ one
observes that they all are derivations with respect to $\mu$ but
$\degs$ and $\degnu$ fail to be derivations with respect to
$\circ$. Instead $\Deg := \degs + 2\degnu$ is a derivation of
$\circ$ and hence $(\WL, \circ)$ is formally $\Deg$-graded and the
corresponding degree is referred to as the total degree. Sometimes
we write ${\mbox{$\mathcal W_k \!
\otimes \! \Lambda$}}$ to denote the elements of total degree
$\geq k$.

In local coordinates we define the differential $\delta : = (1
\otimes \d x^i) i_s (\partial_i)$ which satisfies $\delta^2=0$ and
is a super-derivation of $\circ$. Moreover, there is a homotopy
operator $\delta^{-1}$ satisfying $\delta \delta^{-1} +
\delta^{-1}\delta + \sigma = \id$ where $\sigma: \WL
\to \Cinf{M}[[\nu]]$ denotes the projection onto the part of
symmetric and anti-symmetric degree $0$ and $\delta^{-1} a :=
\frac{1}{k+l} (\d x^i\otimes 1)i_a(\partial_i)a$ for $\degs a = k
a$, $\dega a = l a$ with $k+l \neq 0$ and $\delta^{-1}a:=0$ else.
From a torsion free symplectic connection $\nabla$ on $M$ we obtain
a derivation $\nabla:= (1 \otimes \d x^i) \nabla_{\partial_i}$ of
$\circ$ that satisfies the following identities:
$[\delta,\nabla]=0$, $\nabla^2=
-\frac{1}{\nu}\ad(R)$, where $R:=\frac{1}{4} \omega_{it}R^t_{jkl}
\d x^i \vee \d x^j \otimes \d x^k \wedge \d x^l \in \WL^2$ involves the
curvature of the connection. Moreover we have $\delta R = 0=\nabla
R$ by the Bianchi identities.

Now remember the following facts which are just restatements of
Fedosov's original theorems in \cite[Thm. 3.2, 3.3]{Fed94} resp.
\cite[Thm. 5.3.3]{Fed96}:

For all $\Omega\in\nu Z^2_{\mbox{\rm\tiny dR}}(M)[[\nu]]$ and all
$s\in\W_3$ with $\sigma(s)=0$ there exists a unique element $r\in
\mbox{$\mathcal W_2 \! \otimes \! \Lambda$}^1$ such that
\begin{equation} \label{Genr}
    \delta r = \nabla r - \frac{1}{\nu} r \circ r + R + 1 \otimes
    \Omega \quad \mbox { and } \quad
    \delta^{-1} r = s.
\end{equation}
Moreover $r$ satisfies the formula
\begin{equation} \label{GenrRecus}
    r = \delta s + \delta^{-1} \left( \nabla r - \frac{1}{\nu}
    r \circ r + R + 1 \otimes \Omega\right)
\end{equation}
from which $r$ can be determined recursively. In this case the
Fedosov derivation
\begin{equation} \label{GenFedDerivDef}
    \mathcal D := - \delta + \nabla - \frac{1}{\nu} \ad(r)
\end{equation}
is a super-derivation of anti-symmetric degree $1$ and has square
zero: $\mathcal D^2 = 0$. Furthermore observe that the $\mathcal
D$-cohomology on elements $a$ with positive anti-symmetric degree
is trivial since one has the following homotopy formula $\mathcal D
\mathcal D^{-1} a + \mathcal D^{-1} \mathcal Da =a$, where
$\mathcal D^{-1} a := -\delta^{-1}\left(\frac{1}{\id- [\delta^{-1},
\nabla -\frac{1}{\nu}\ad(r)]} a \right)$
(cf. \cite[Thm. 5.2.5]{Fed96}).

Then for any $f \in \Cinf{M}[[\nu]]$ there exists a unique element
$\tau(f) \in \ker (\mathcal D) \cap \W$ such that $\sigma (\tau(f))
= f$ and $\tau: \Cinf{M}[[\nu]] \to \ker (\mathcal D)
\cap\W$ is $\mathbb C[[\nu]]$-linear and referred to as the
Fedosov-Taylor series corresponding to $\mathcal D$. In addition
$\tau(f)$ can be obtained recursively for $f\in \Cinf{M}$ from
      \begin{equation} \label{GenTaylorRecurs}
      \tau(f)= f +\delta^{-1}\left(\nabla \tau(f)-\frac{1}{\nu}
      \ad(r)\tau(f)
      \right).
      \end{equation}
Using $\mathcal D^{-1}$ one can also write $\tau(f) = f - \mathcal
D^{-1} (1 \otimes \d f)$. Since $\mathcal D$ as constructed above
is a $\circ$-super-derivation $\ker(\mathcal D) \cap \W$ is a
$\circ$-sub-algebra and a new associative product $*$ for
$\Cinf{M}[[\nu]]$, which turns out to be a star product, is defined
by pull-back of $\circ$ via $\tau$.

Observe that in (\ref{Genr}) we allowed for an arbitrary element $s
\in \W$ with $\sigma(s)=0$ that contains no terms of total degree
lower than $3$, as normalization condition for $r$, i.e.
$\delta^{-1}r=s$ instead of the usual equation $\delta^{-1}r=0$. In
the following we shall refer to the associative product $*$ defined
above as the Fedosov star product (corresponding to
$(\nabla,\Omega,s)$).

Now we shall give a very convenient description of all derivations
of the star product $*$ that will prove very useful for our further
considerations. To this end we consider appropriate fibrewise
quasi-inner derivations of the shape
\begin{equation}
\mathrm D_h = -\frac{1}{\nu} \ad (h),
\end{equation}
where $h \in \W$ and without loss of generality we assume
$\sigma(h)=0$. Our aim is to define $\mathbb C[[\nu]]$-linear
derivations of $*$ by $\Cinf{M}[[\nu]] \ni f \mapsto
\sigma(\mathrm D_h
\tau(f))$ but for an arbitrary element $h\in \W$ with $\sigma(h)=
0$ this mapping fails to be a derivation as $\mathrm D_h$ does not
map elements of $\ker(\mathcal D)\cap \W$ to elements of
$\ker(\mathcal D)\cap \W$. In order to achieve this one must have
that $\mathcal D$ and $\mathrm D_h$ super-commute. As $\mathcal D$
is a $\mathbb C[[\nu]]$-linear $\circ$-super-derivation we
obviously have
\[
[\mathcal D, \mathrm D_h] = - \frac{1}{\nu} \ad (\mathcal D h)
\]
and hence obviously $\mathcal D h$ must be central, i.e. $\mathcal
D h$ has to be of the shape $1 \otimes A$ with $A \in
\Ginf{T^*M}[[\nu]]$ to have $[\mathcal D, \mathrm D_h]=0$. From
$\mathcal D^2 =0$ we get that the necessary condition for the
solvability of the equation $\mathcal D h = 1 \otimes A$ is the
closedness of $A$ since $\mathcal D (1 \otimes A ) = 1 \otimes
\d A$. But as the $\mathcal D$-cohomology is trivial on elements
with positive anti-symmetric degree this condition is also
sufficient for the solvability of the equation $\mathcal D h = 1
\otimes A$ and we get the following statement.

\begin{LEMMA}\label{DerLem}
\begin{enumerate}
\item
For all formal series $A\in \Ginf{T^*M}[[\nu]]$ of closed one-forms
on $M$ there is a uniquely determined element $h_A \in \W$ such
that $\mathcal D h_A = 1 \otimes A$ and $\sigma (h_A) =0$. Moreover
$h_A$ is explicitly given by
\begin{equation}
h_A = \mathcal D^{-1} (1 \otimes A).
\end{equation}
\item
For all $A \in Z^1_{\mbox{\rm\tiny dR}}(M)[[\nu]]$ the mapping
$\mathsf D_A : \Cinf{M}[[\nu]] \to \Cinf{M}[[\nu]]$, where
\begin{equation}
\mathsf D_A f := \sigma (\mathrm D_{h_A}\tau(f)) = \sigma\left(
-\frac{1}{\nu}\ad(h_A)\tau(f) \right)
\end{equation}
for $f \in \Cinf{M}[[\nu]]$ defines a $\mathbb C[[\nu]]$-linear
derivation of $*$ and hence this construction yields a mapping
$Z^1_{\mbox{\rm\tiny dR}}(M)[[\nu]] \ni A \mapsto \mathsf D_A \in
\mathrm{Der}_{\mathbb C[[\nu]]}(\Cinf{M}[[\nu]],*)$.
\end{enumerate}
\end{LEMMA}
\begin{PROOF}
The fact that $h_A = \mathcal D^{-1}(1 \otimes A)$ satisfies
$\mathcal D h_A = 1 \otimes A$ is obvious from the homotopy formula
for $\mathcal D$ and the closedness of $A$. In addition we have
$\sigma(h_A)=0$ since $\mathcal D^{-1}$ raises the symmetric degree
at least by $1$. For the uniqueness of $h_A$ let $\widetilde{h_A}$
be another solution of the equations above, then we obviously have
$\mathcal D(h_A - \widetilde{h_A})=0$ and hence $h_A -
\widetilde{h_A} = \tau(\varphi)$ for some $\varphi \in
\Cinf{M}[[\nu]]$. Applying $\sigma$ to this equation one gets
$\varphi=0$, since $\sigma(h_A)=\sigma(\widetilde{h_A})=0$ and
$\sigma(\tau(\varphi))=\varphi$, and hence $h_A
= \widetilde{h_A}$ proving that $h_A$ is uniquely determined by the
above equations. For the proof of ii.) we just observe that the
equation $[\mathcal D, \mathrm D_{h_A}]=0$ which is fulfilled
according to i.) implies that $\mathrm D_{h_A} \tau (f) =
\tau(\mathsf D_A f)$ for all $f \in \Cinf{M}[[\nu]]$. Using this
equation and the obvious fact that $\mathrm D_{h_A}$ is a
derivation of $\circ$ it is straightforward to see using the very
definition of $*$ that $\mathsf D_A$ as defined above is a
derivation of $*$. The $\mathbb C[[\nu]]$-linearity of $\mathsf
D_A$ is also evident from the $\mathbb C[[\nu]]$-linearity of
$\tau$.
\end{PROOF}

Furthermore we now are in the position to show that one even
obtains all $\mathbb C[[\nu]]$-linear derivations of $*$ by varying
$A$ in the derivations $\mathsf D_A$ constructed above.

\begin{PROPOSITION}\label{DerBijProp}
The mapping
\[
Z^1_{\mbox{\rm\tiny dR}}(M)[[\nu]] \ni A \mapsto \mathsf D_A \in
\mathrm{Der}_{\mathbb C[[\nu]]}(\Cinf{M}[[\nu]],*)
\]
defined in Lemma \ref{DerLem} is a bijection. Moreover, $\mathsf
D_{\d f}$ is a quasi-inner derivation for all $f \in
\Cinf{M}[[\nu]]$, i.e. $\mathsf D_{\d f}= \frac{1}{\nu} \ad_* (f)$
and the induced mapping $[A] \mapsto [\mathsf D_A]$ from
$H^1_{\mbox{\rm\tiny dR}}(M)[[\nu]]\cong Z^1_{\mbox{\rm\tiny
dR}}(M)[[\nu]]/B^1_{\mbox{\rm\tiny dR}}(M)[[\nu]]$ to
$\mathrm{Der}_{\mathbb C[[\nu]]}(\Cinf{M}[[\nu]],*)/
\mathrm{Der}^{\mbox{\rm\tiny qi}}_{\mathbb C[[\nu]]}
(\Cinf{M}[[\nu]],*)$ the space of $\mathbb C[[\nu]]$-lin\-e\-ar
derivations of $*$ modulo the quasi-inner derivations, also is
bijective.
\end{PROPOSITION}
\begin{PROOF}
First we prove the injectivity of the mapping $A\mapsto \mathsf
D_A$. To this end let $\mathsf D_A = \mathsf D_{A'}$ then we get
from $\mathrm D_{h_A} \tau(f) = \tau (\mathsf D_A f)$ and from the
analogous equation for $A'$ that $\ad(h_A - h_{A'})\tau(f)=0$ for
all $f\in \Cinf{M}[[\nu]]$ and hence $h_A - h_{A'}$ must be central
(since it commutes with all Fedosov-Taylor series), i.e. we have
$h_A - h_{A'} = g_{A,A'}\in \Cinf{M}[[\nu]]$. But with
$\sigma(h_A)=\sigma(h_{A'})=0$ this implies $g_{A,A'}=0$ and hence
$h_A = h_{A'}$ such that we get $1\otimes A = \mathcal D h_A
= \mathcal D h_{A'} =1\otimes A'$ proving the injectivity. For the
surjectivity we start with an arbitrary derivation $\mathsf D$ of
$*$ and want to find closed one-forms $A_i$ such that $\mathsf D =
\sum_{i=0}^\infty \nu^i \mathsf D_{A_i}$ inductively. Assume that we have found
such one-forms for $0 \leq i \leq k-1$ such that $\mathsf D' =
\mathsf D - \sum_{i=0}^{k-1} \nu^i \mathsf D_{A_i}$ which obviously
is again a derivation of $*$ is of the shape $\mathsf D'=
\sum_{i=k}^\infty \nu^i \mathsf D'_i$. The $k$th order in $\nu$ of
the equation $\mathsf D'(f * g) = (\mathsf D'f) * g + f
* (\mathsf D' g)$ for $f,g \in \Cinf{M}$ yields that $\mathsf D'_k$
is a vector field $X_k \in \Ginf{TM}$. Considering the
anti-symmetric part of $\mathsf D'(f * g) = (\mathsf D'f) * g + f *
(\mathsf D' g)$ at order $k+1$ of $\nu$ we get that this vector
field is symplectic, i.e. $\Lie_{X_k} \omega =0$ and because of the
Cartan formula $A_k
:= - i_{X_k} \omega$ defines a closed one-form on $M$. Considering
the derivation $\mathsf D_{A_k}$ it is a straightforward
computation using the explicit construction above to show that
$\mathsf D_{A_k} f = X_k (f) + O(\nu)$ for all $f \in \Cinf{M}$.
But then $\mathsf D' - \nu^k \mathsf D_{A_k}$ is again a derivation
of $*$ that starts in order $k+1$ of $\nu$ and hence the
surjectivity follows by induction. The fact that $\mathsf D_{\d f}
= \frac{1}{\nu} \ad_* (f)$ for all $f\in \Cinf{M}[[\nu]]$ is
obvious from the observation that $\tau(f) = f - \mathcal D^{-1} (1
\otimes \d f)$ and the obvious fact that $\ad (f)=0$. From the
above, the well-definedness of the mapping $[A]\mapsto [\mathsf
D_A]$ follows and the bijectivity is a direct consequence of the
bijectivity of the mapping $A \mapsto \mathsf D_A$.
\end{PROOF}

\begin{REMARK}
Actually it is well-known that for an arbitrary star product
$\star$ on a symplectic manifold the space of $\mathbb
C[[\nu]]$-linear derivations is in bijection with
$Z^1_{\mbox{\rm\tiny dR}}(M)[[\nu]]$ and that the quotient space of
these derivations modulo the quasi-inner derivations is in
bijection with $H^1_{\mbox{\rm\tiny dR}}(M)[[\nu]]$ (cf. \cite[Thm.
4.2]{BerCahGut97}, observe that the proof given above is just an
adaption of the idea of the general proof to our special situation)
but the remarkable thing about Fedosov star products is that these
bijections can be explicitly expressed in terms of $\mathcal D$
resp. $\mathcal D^{-1}$ in a very lucid way which will be useful in
the following.
\end{REMARK}

To conclude this section we shall remove some redundancy in the
description of the star products $*$ by $(\nabla,\Omega,s)$. This
will ease the more detailed analysis in the following section. To
this end we shall recall some well-known facts about symplectic
torsion free connections on $(M,\omega)$. Given two such
connections say $\nabla$ and $\nabla'$ it is obvious that
$S^{\nabla-\nabla'}(X,Y):= \nabla_X Y - \nabla'_X Y$ where $X,Y\in
\Ginf{TM}$ defines a symmetric tensor field $S^{\nabla - \nabla'}
\in \Ginf{\bigvee^2 T^*M \otimes TM}$ on $M$. Defining
$\sigma^{\nabla-\nabla'}(X,Y,Z):= \omega(S^{\nabla -
\nabla'}(X,Y),Z)$ it is easy to see that $\sigma^{\nabla-\nabla'}
\in \Ginf{\bigvee^3 T^*M}$ is a totally symmetric tensor field.
Vice versa given an arbitrary element $\sigma\in\Ginf{\bigvee^3
T^*M}$ and a symplectic torsion free connection $\nabla$ and
defining $S^\sigma\in \Ginf{\bigvee^2 T^*M \otimes TM}$ by
$\sigma(X,Y,Z) = \omega(S^\sigma(X,Y),Z)$ then $\nabla^\sigma$
defined by $\nabla^\sigma_X Y := \nabla_X Y - S^\sigma(X,Y)$ again
is a symplectic torsion free connection and all such connections
can be obtained this way by varying $\sigma$. Using these relations
we shall compare the corresponding mappings $\nabla$ and $\nabla'$
on $\WL$ in the following lemma.

\begin{LEMMA}
With the notations from above we have
\begin{equation}
\label{NabNabStrDifEq}
\nabla- \nabla' = - (\d x^j \otimes \d x^i) i_s(S^{\nabla-\nabla'}
(\partial_i,\partial_j)) = \frac{1}{\nu} \ad (T^{\nabla-\nabla'}),
\end{equation}
where $T^{\nabla-\nabla'}\in \Ginf{\mbox{$\bigvee$}^2T^*M \otimes
T^*M}\subseteq\WL^1$ is defined by
$T^{\nabla-\nabla'}(Z,Y;X):=\sigma^{\nabla-\nabla'}(X,Y,Z)=
\omega(S^{\nabla-\nabla'}(X,Y),Z)$. Moreover $T^{\nabla-\nabla'}$
satisfies the equations
\begin{equation}\label{TNabNabStrEq}
\delta T^{\nabla-\nabla'}= 0 \qquad \textrm{ and } \qquad
\begin{array}{c} \nabla T^{\nabla-\nabla'} = R' -R +
\frac{1}{\nu} T^{\nabla-\nabla'}\circ T^{\nabla-\nabla'}\\ \nabla'
T^{\nabla-\nabla'} = R' -R - \frac{1}{\nu} T^{\nabla-\nabla'}\circ
T^{\nabla-\nabla'},
\end{array}
\end{equation}
where $R = \frac{1}{4} \omega_{it} R^t_{jkl}\d x^i\vee \d
x^j\otimes \d x^k\wedge \d x^l$ and $R'=\frac{1}{4} \omega_{it}
{R'}^t_{jkl}\d x^i\vee \d x^j\otimes \d x^k\wedge \d x^l$ denote
the corresponding elements of $\WL^2$ that are built from the
curvature tensors of $\nabla$ and $\nabla'$.
\end{LEMMA}
\begin{PROOF}
The proof of (\ref{NabNabStrDifEq}) is a straightforward
computation using the very definitions from above. The first
identity in (\ref{TNabNabStrEq}) directly follows from
(\ref{NabNabStrDifEq}) and $[\delta,\nabla]=[\delta,\nabla']=0$.
The other identities in (\ref{TNabNabStrEq}) are also easily
obtained squaring equation (\ref{NabNabStrDifEq}).
\end{PROOF}

Now we are in the position to compare two Fedosov derivations
$\mathcal D$ and $\mathcal D'$ resp. the induced star products $*$
and $*'$ obtained from $(\nabla, \Omega, s)$ and $(\nabla',\Omega',
s')$.

\begin{PROPOSITION}
The Fedosov derivations $\mathcal D$ and $\mathcal D'$ coincide if
and only if $T^{\nabla-\nabla'} - r + r' = 1 \otimes \vartheta$
where $\vartheta \in \nu \Ginf{T^*M}[[\nu]]$ which is equivalent to
\begin{equation}\label{DD'equilEq}
\sigma^{\nabla -\nabla'} \otimes 1 - s + s' = \vartheta \otimes 1
\quad\textrm{ and } \quad \Omega - \Omega' = \d \vartheta.
\end{equation}
\end{PROPOSITION}
\begin{PROOF}
Writing down the definitions of $\mathcal D$ and $\mathcal D'$
using equation (\ref{NabNabStrDifEq}) the first equivalence is
obvious since $T^{\nabla - \nabla'} -r +r'$ is central in
$(\WL,\circ)$ if and only if $\mathcal D = \mathcal D'$. For the
proof of the second equivalence first assume that we have
$T^{\nabla-\nabla'} - r + r' = 1 \otimes \vartheta$. Applying
$\delta^{-1}$ to this equation and using the normalization
condition on $r$ and $r'$ we obtain the first equation in
(\ref{DD'equilEq}) since $\delta^{-1}T^{\nabla -\nabla'} =
\sigma^{\nabla-\nabla'}\otimes 1$. In order to obtain the
second equation in (\ref{DD'equilEq}) we apply $\delta$ to
$T^{\nabla-\nabla'} - r + r' = 1 \otimes \vartheta$ and a
straightforward computation using the equations for $r$ and $r'$
together with the identities from (\ref{TNabNabStrEq}) yields the
stated result. To prove that the converse is also true assume that
the equations in (\ref{DD'equilEq}) are satisfied and define $B:=
r-r' - T^{\nabla-\nabla'} + 1 \otimes \vartheta \in
{\mbox{$\mathcal W_2 \! \otimes \! \Lambda$}}^1$. Then again a
straightforward computation yields that $B$ satisfies $\mathcal D B
= - \frac{1}{\nu} B \circ B$ and $\delta^{-1} B=0$ such that the
homotopy formula for $\delta$ together with $\sigma(B)=0$ implies
that $B$ is the unique fixed point of the mapping ${\mbox{$\mathcal
W_2 \! \otimes \! \Lambda$}}^1 \ni a \mapsto
\delta^{-1}\left(\nabla a - \frac{1}{\nu} \ad (r) a + \frac{1}{\nu}
a\circ a \right) \in {\mbox{$\mathcal W_2 \! \otimes \!
\Lambda$}}^1$. But $0$ trivially is a fixed point of this mapping
and hence uniqueness implies that $B=0$ proving the other direction
of the second stated equivalence.
\end{PROOF}

As an important direct consequence of this proposition we get:
\begin{DEDUCTION}\label{RedunDed}
For every Fedosov star product $*$ obtained from $(\nabla,\Omega,
s)$ with $s \in \W_3$ there is a connection $\nabla'$, a formal
series $\Omega'$ of closed two-forms and an element $s' \in \W_4$
without terms of symmetric degree $1$ such that the star product
obtained from $(\nabla',\Omega',s')$ coincides with $*$, and hence
we may without loss of generality restrict to such normalization
conditions when varying the connection and the formal series of
closed two-forms arbitrarily.
\end{DEDUCTION}
\begin{PROOF}
We write $s= s' + \sigma \otimes 1 - \vartheta \otimes 1$ and the
preceding proposition states that $\mathcal D$ coincides with
$\mathcal D'$ (and hence the corresponding star products coincide)
where $\mathcal D'$ is obtained from $\Omega'
= \Omega - \d\vartheta$ and $\nabla' = \nabla - \frac{1}{\nu}\ad
(\delta(\sigma\otimes 1))$.
\end{PROOF}

\section{Symplectic Vector Fields as Derivations of $*$}
\label{SymVecDerSec}
Throughout this and the following section let $*$ denote the
Fedosov star product obtained from $(\nabla,\Omega,s)$ as in
Section \ref{PrelimSec} where in view of Deduction \ref{RedunDed}
we may assume that $s\in \W_4$ contains no part of symmetric degree
$1$. Furthermore $X\in \Ginf{TM}$ shall always denote a symplectic
vector field on $(M,\omega)$ and the space of all these vector
fields shall be denoted by $\symp:=\{Y \in \Ginf{TM}\,|\, \Lie_Y
\omega =0\}$. It seems to be folklore and actually is not very hard to prove
that the conditions $[\Lie_X, \nabla]=0$, $\Lie_X \Omega =0=\Lie_X
s$ are sufficient to guarantee that the Lie derivative with respect
to $X$ is a derivation of $*$. Besides providing a very simple
proof of this fact, our aim in this section is to prove that the
converse is also true, i.e. the conditions given above are also
necessary to have that $X$ defines a derivation of $*$. Moreover,
we find an additional cohomological condition involving $\omega$,
$\Omega$ and $X$ that is equivalent to $\Lie_X$ being even a
quasi-inner derivation.

As an important tool we need the deformed Cartan formula (cf.
\cite[Appx. A]{Neu99}) that relates the Lie derivative with respect
to a symplectic vector field $X$ with the Fedosov derivation
$\mathcal D$.

\begin{LEMMA}
For all $X\in \symp$ the Lie derivative $\Lie_X$ can be expressed
in the following manner:
\begin{equation}\label{CarForEq}
\Lie_X = \mathcal D i_a(X) + i_a(X)\mathcal D - \frac{1}{\nu}
\ad\left(\theta_X \otimes 1 + \frac{1}{2}D \theta_X \otimes 1 -
i_a(X)r \right),
\end{equation}
where $D:= \d x^i \vee \nabla_{\partial_i}$ denotes the operator of
symmetric covariant derivation and the closed one-form $\theta_X$
is defined by $\theta_X:= i_X \omega$.
\end{LEMMA}
\begin{PROOF}
Since the Lie derivative is a local operator it suffices to prove
the above identity over any contractible open subset $U$ of $M$.
But as $X$ is symplectic it is locally Hamiltonian, i.e. over $U$
there is a function $f \in \Cinf{U}$ such that $X|_U = X_f$ resp.
$\d f = \theta_X |_U$. For Hamiltonian vector fields the Cartan
formula as above was proved in \cite[Prop. 5]{Neu99} and hence
equation (\ref{CarForEq}) is valid for all symplectic vector fields
$X\in\symp$.
\end{PROOF}

As an immediate consequence of the preceding lemma we have:
\begin{LEMMA}
For $X\in \symp$ the Lie derivative $\Lie_X$ is a derivation with
respect to $\circ$. In addition we have $[\delta,\Lie_X]=
[\delta^{-1},\Lie_X]=0$.
\end{LEMMA}
\begin{PROOF}
The first statement of the lemma is obvious from equation
(\ref{CarForEq}) and the commutation relations follow from the fact
that $\Lie_X$ is compatible with contractions and preserves the
symmetric and the anti-symmetric degree.
\end{PROOF}

After these rather technical preparations we get:
\begin{PROPOSITION}\label{DerProp}
Let $X\in \symp$ then $\Lie_X$ is a derivation of $*$ if and only
if $[\Lie_X,\mathcal D]=0$ which is equivalent to the existence of
a formal series $A_X\in
\Ginf{T^*M}[[\nu]]$ of closed one-forms such that $\mathcal D
\left(\theta_X\otimes 1 + \frac{1}{2}D \theta_X \otimes 1 -i_a(X)r
\right)= 1 \otimes A_X$.
\end{PROPOSITION}
\begin{PROOF}
First let us assume that $[\Lie_X,\mathcal D]=0$ then the obvious
equation $\Lie_X \circ \sigma = \sigma \circ \Lie_X$ implies that
$\Lie_X \tau(f) = \tau (\Lie_X f)$ for all $f\in
\Cinf{M}[[\nu]]$. But with this equation and the fact that $\Lie_X$
is a derivation of $\circ$ it is straightforward to prove that
$\Lie_X$ is a derivation of $*$. Assuming that $\Lie_X$ is a
derivation of $*$ Proposition \ref{DerBijProp} implies that there
is a formal series $A_X$ of closed one-forms on $M$ such that
$\Lie_X f = \sigma\left( - \frac{1}{\nu} \ad (\mathcal
D^{-1}(1\otimes A_X))\tau(f) \right)$ but on the other hand the
deformed Cartan formula yields $\Lie_X f = \sigma\left(
-\frac{1}{\nu}\ad\left(\theta_X \otimes 1+ \frac{1}{2}D\theta_X
\otimes 1 - i_a(X)r \right)\tau(f)\right)$ and hence $\mathcal
D^{-1}(1\otimes A_X)-\left(\theta_X \otimes 1+ \frac{1}{2}D\theta_X
\otimes 1 - i_a(X)r \right)$ has to be central, i.e. a formal
function. Observing that $\mathcal D^{-1}$ raises the symmetric
degree at least by $1$ and that $r$ contains no part of symmetric
degree $0$ which is due to the special shape of the normalization
condition this implies $\mathcal D^{-1}(1\otimes A_X)=
\left(\theta_X \otimes 1+ \frac{1}{2}D\theta_X
\otimes 1 - i_a(X)r \right)$. Applying $\mathcal D$ to this
equation and using the homotopy formula for $\mathcal D$ together
with the fact that $A_X$ is closed we get $\mathcal D
\left(\theta_X\otimes 1 + \frac{1}{2}D \theta_X \otimes 1 -i_a(X)r
\right)= 1 \otimes A_X$. Assuming finally that this equation is
fulfilled, the deformed Cartan formula together with $\mathcal
D^2=0$ obviously implies $[\Lie_X,\mathcal D]=0$ since $1\otimes
A_X$ is central and hence the proposition is proved.
\end{PROOF}

We shall now go on by analysing the condition
\begin{equation}\label{DerCondEq}
\mathcal D
\left(\theta_X\otimes 1 + \frac{1}{2}D \theta_X \otimes 1 -i_a(X)r
\right)= 1 \otimes A_X,\quad\textrm{ where } \quad \d A_X=0
\end{equation}
in more detail in order to find out whether it gives rise to
conditions on $(\nabla,\Omega,s)$ and $X$.

\begin{LEMMA}\label{CondLem1}
For all symplectic vector fields $X\in \symp$ we have
\begin{equation}\label{CondEq1}
\mathcal D
\left(\theta_X\otimes 1 + \frac{1}{2}D \theta_X \otimes 1 -i_a(X)r
\right)=- 1 \otimes \theta_X + \nabla\left(\frac{1}{2}D \theta_X
\otimes 1 \right) - \Lie_X r - i_a(X)R - 1 \otimes i_X \Omega.
\end{equation}
\end{LEMMA}
\begin{PROOF}
The proof of this equation is a straightforward computation using
the equation that is solved by $r$ and the deformed Cartan formula
(\ref{CarForEq}) once again.
\end{PROOF}

Next we shall need some detailed formulas that describe
$[\nabla,\Lie_X]$ in order to simplify the result of the above
Lemma. The proofs of the following two lemmas are just slight
variations of the proofs of \cite[Lemma 3 and Lemma 4]{Neu99}.

\begin{LEMMA}\label{SXLem}
For all $X\in \symp$ the mapping $[\nabla,\Lie_X]$ enjoys the
following properties:
\begin{enumerate}
\item
In local coordinates one has
\begin{equation} \label{nablaLiecomEq}
    [\nabla,\Lie_X] = (\d x^j \otimes \d x^i)
    i_s((\Lie_X \nabla)_{\partial_i}\partial_j)= (\d x^j
    \otimes \d x^i) i_s(S_X(\partial_i, \partial_j)),
\end{equation}
where the tensor field $S_X \in
\Ginf{T^*M \otimes T^*M \otimes TM}$ is defined by
\begin{equation} \label{SXdef}
    S_X(\partial_i,\partial_j)=(\Lie_X
    \nabla)_{\partial_i} \partial_j:= \Lie_X
    \nabla_{\partial_i}\partial_j - \nabla_{\partial_i}
    \Lie_X\partial_j - \nabla_{\Lie_X\partial_i}\partial_j =
    R(X, \partial_i) \partial_j +
    \nabla^{(2)}_{(\partial_i,\partial_j)}X.
\end{equation}
\item
$S_X$ as defined above is symmetric, i.e. $S_X \in
\Ginf{\bigvee^2 T^*M \otimes TM}$.
\item
For all $U, V, W \in \Ginf{T M}$ we have $\omega(W, S_X(U,V))= -
\omega(S_X(U,W),V)$.
\end{enumerate}
\end{LEMMA}
Now the tensor field $S_X$ naturally gives rise to an element
$T_X\in \Ginf{\bigvee^2T^*M\otimes T^*M}$ of $\WL^1$ of symmetric
degree $2$ and anti-symmetric degree $1$ by
\begin{equation}\label{TXDefEq}
T_X(W,U; V) := \omega (W, S_X(V,U))
\end{equation}
and we have:
\begin{LEMMA}\label{TXLem}
The tensor field $T_X$ as defined in (\ref{TXDefEq}) satisfies the
following equations:
\begin{enumerate}
\item
    $\frac{1}{\nu}\ad (T_X) = [\nabla,\Lie_X]$,
\item
    $T_X = i_a(X) R - \nabla\left(\frac{1}{2}D
    \theta_X \otimes 1\right)$,
\item
    $\delta T_X = 0 \quad\textrm{ and } \quad\nabla T_X =
    \Lie_X R$.
\end{enumerate}
\end{LEMMA}

From the preceding lemma we find that the result of Lemma
\ref{CondLem1} simplifies to
\begin{equation}\label{CondEq2}
\mathcal D
\left(\theta_X\otimes 1 + \frac{1}{2}D \theta_X \otimes 1 -i_a(X)r
\right)=- 1 \otimes \theta_X - T_X - \Lie_X r - 1 \otimes i_X
\Omega.
\end{equation}
Finally we have to find equations that determine $\Lie_X r$ in
order to analyse equation (\ref{DerCondEq}).

\begin{LEMMA}
Let $X$ denote a symplectic vector field then $\Lie_X r$ satisfies
the equations
\begin{equation}\label{LieXrEq}
\delta \Lie_X r =  \nabla\Lie_X r
-\frac{1}{\nu} \ad(r)\Lie_X r - \frac{1}{\nu}\ad(T_X) r
+ \Lie_X R + 1 \otimes \d i_X \Omega
\textrm{ and } \delta^{-1}\Lie_X r = \Lie_X s
\end{equation}
from which $\Lie_X r$ is uniquely determined and can be computed
recursively from
\[
\Lie_X r = \delta \Lie_X s + \delta^{-1}\left(
\nabla\Lie_X r - \frac{1}{\nu}\ad(r) \Lie_X r - \frac{1}{\nu}
\ad(T_X)r + \Lie_X R + 1 \otimes \d i_X \Omega
\right).
\]
\end{LEMMA}
\begin{PROOF}
For the proof of (\ref{LieXrEq}) one just has to apply $\Lie_X$ to
the equations that determine $r$ and to use the commutation
relations of the involved mappings. From these equations it is
straightforward to find the recursion formula for $\Lie_X$ using
the homotopy formula for $\delta$. Using statement iii.) of Lemma
\ref{TXLem} the argument for the uniqueness of the solution of
these equations is completely analogous to the one used to prove
the uniqueness of $r$ and hence we leave it to the reader.
\end{PROOF}

After all these preparations we are in the position to formulate
the main results of this section.

\begin{THEOREM}\label{DerTheo}
Let $X$ be a symplectic vector field and let $*$ be the Fedosov
star product corresponding to $(\nabla,\Omega,s)$, where $s\in
\W_4$ contains no part of symmetric degree $1$. Then,
$\Lie_X$ is a derivation of $*$ if and only if $T_X=0$, $\Lie_X
\Omega =0$ and $\Lie_X s=0$, i.e. if and only if $X$ is affine with
respect to $\nabla$ and $s$ and $\Omega$ are invariant with respect
to $X$.
\end{THEOREM}
\begin{PROOF}
First let $T_X=0=\Lie_X \Omega = \Lie_X s$ then we have $\Lie_X R
=\nabla T_X =0$ and $\d i_X \Omega=0$ and hence $\Lie_X r =
\delta^{-1}\left( \nabla \Lie_X r - \frac{1}{\nu} \ad(r) \Lie_X r
\right)$. But this implies $\Lie_X r=0$ and then obviously
$[\mathcal D,\Lie_X]=\frac{1}{\nu}\ad (T_X + \Lie_X r)=0$ such that
Proposition \ref{DerProp} implies that $\Lie_X$ is a derivation of
$*$. To prove the converse we again use Proposition \ref{DerProp}
which says that in case $\Lie_X$ is a derivation of $*$ there is a
formal series $A_X$ of closed one-forms on $M$ such that $\mathcal
D
\left(\theta_X\otimes 1 + \frac{1}{2}D \theta_X \otimes 1 -i_a(X)r
\right)=1 \otimes A_X$. Together with equation (\ref{CondEq2}) this
yields $\Lie_X r = - ( 1\otimes (\theta_X + A_X + i_X\Omega)
+T_X)$. Applying $\delta^{-1}$ to this equation and using the
second equation in (\ref{LieXrEq}) we get
\[
\Lie_X s = - (\theta_X + A_X + i_X\Omega)\otimes 1 - \delta^{-1} T_X.
\]
Now $s$ and hence $\Lie_X s$ is in $\W_4$ and has no part of
symmetric degree $1$ such that this equation implies $\Lie_X s=0$,
$\theta_X + A_X + i_X\Omega=0$ and $\delta^{-1}T_X=0$. Since
$\theta_X$ and $A_X$ are closed the second of these equations
implies $0=\d i_X\Omega=\Lie_X \Omega$ and using the homotopy
formula for $\delta$ together with $\delta T_X=0$ the last equation
yields $T_X=0$ which is equivalent to $X$ being affine with respect
to $\nabla$ according to the Lemmas \ref{SXLem} and \ref{TXLem}.
Finally one can insert the above expression for $\Lie_X r$ into the
first equation in (\ref{LieXrEq}) which turns out to be satisfied
identically, which is just a check for consistency.
\end{PROOF}

Finally we can give an additional condition for $\Lie_X$ to be even
a quasi-inner derivation of $*$ which is originally due to Gutt
\cite{Gut02}.

\begin{PROPOSITION}\label{quainnDerProp}
Let $X$ be a symplectic vector field such that $\Lie_X$ is a
derivation of $*$ then $\Lie_X$ is even quasi-inner if and only if
there is a formal function $f \in \Cinf{M}[[\nu]]$ such that
\begin{equation}\label{qinnEq}
\d f = \theta_X + i_X \Omega = i_X (\omega +\Omega)
\end{equation}
and then $\Lie_X =\Lie_{X_{f_0}}= - \frac{1}{\nu}
\ad_*(f)$, where we have written $f = f_0
+ f_+$ with $f_0\in
\Cinf{M}$ and $f_+ \in \nu \Cinf{M}[[\nu]]$.
\end{PROPOSITION}
\begin{PROOF}
From equation (\ref{CarForEq}) it is obvious that $\Lie_X$ is
quasi-inner if and only if there is a formal function $f
\in \Cinf{M}[[\nu]]$ such that $\tau (f) =
f +\theta_X\otimes 1 + \frac{1}{2}D \theta_X \otimes 1
-i_a(X)r$ but using equation (\ref{CondEq2}) together with $T_X=0$,
$\Lie_X r=0$ and $\mathcal D f = 1 \otimes \d f$ this is equivalent
to (\ref{qinnEq}). In fact the necessary condition for the
solvability of this equation is fulfilled since $i_X\Omega$ is
closed according to Theorem \ref{DerTheo} and $\theta_X$ is closed
as $X$ is symplectic. Moreover, observe that the zeroth order in
$\nu$ of (\ref{qinnEq}) just means that $X$ is Hamiltonian with
Hamiltonian function $f_0$ and hence the second statement of the
Proposition is immediate.
\end{PROOF}

\section{$\mathfrak g$-invariant
Star Products $*$ and Quantum Momentum Mappings}
\label{ginvQMMapSec}
In this section we shall use the results of Theorem \ref{DerTheo}
to find necessary and sufficient conditions for the star product
$*$ to be invariant with respect to a Lie algebra action.
Furthermore Proposition \ref{quainnDerProp} gives criteria for the
existence of a quantum Hamiltonian and with some little more effort
we shall find a last condition which is necessary and sufficient
for this quantum Hamiltonian to define a quantum momentum mapping
for $*$.

First let us recall some definitions from \cite{Xu98}. Let us
consider a finite dimensional real or complex Lie algebra
$\mathfrak g$ and let $X_{\cdot} :
\mathfrak g \to \symp : \xi \mapsto X_\xi$ denote a Lie
algebra anti-homomorphism, i.e. $[X_\xi, X_\eta]= - X_{[\xi,\eta]}$
for all $\xi,\eta \in \mathfrak g$. Then obviously $\varrho(\xi)f:=
- \Lie_{X_\xi} f$ defines a Lie algebra action of $\mathfrak g$ on
$\Cinf{M}$ that naturally extends to a Lie algebra action on
$\Cinf{M}[[\nu]]$.

\begin{DEFINITION}
With the notations from above a star product $\star$ is called
$\mathfrak g$-invariant in case $\varrho(\xi)$ is a derivation of
$\star$ for all $\xi \in \mathfrak g$.
\end{DEFINITION}

From Theorem \ref{DerTheo} we obviously get:
\begin{DEDUCTION}\label{ginvDed}
The Fedosov star product $*$ constructed from $(\nabla,\Omega,s)$,
where $s\in \W_4$ contains no part of symmetric degree $1$, is
$\mathfrak g$-invariant if and only if $X_\xi$ is affine with
respect to $\nabla$ for all $\xi \in \mathfrak g$, i.e.
$[\nabla,\Lie_{X_\xi}]=0\,\forall \xi\in \mathfrak g$ and $\Omega$
and $s$ are invariant with respect to $X_\xi$ for all $\xi
\in \mathfrak g$, i.e. $\d i_{X_\xi} \Omega =
\Lie_{X_\xi}\Omega = 0 = \Lie_{X_\xi}s\,\forall
\xi\in \mathfrak g$.
\end{DEDUCTION}

Let us introduce some notation: Considering some complex vector
space $V$ endowed with a representation $\pi: \mathfrak g \to
\mathsf{Hom}(V,V)$ of the Lie algebra $\mathfrak g$ in $V$ we
denote the space of $V$-valued $k$-multilinear alternating forms on
$\mathfrak g$ by $C^k(\mathfrak g,V)$ and the corresponding
Chevalley-Eilenberg differential shall be denoted by $\delta_\pi :
C^\bullet(\mathfrak g,V)\to C^{\bullet +1}(\mathfrak g,V)$.
Moreover the spaces of the corresponding cocycles and coboundaries
resp. the corresponding cohomology spaces shall be denoted by
$Z^k_\pi(\mathfrak g,V)$ and $B^k_\pi(\mathfrak g,V)$ resp.
$H^k_\pi(\mathfrak g,V)$.

Now the Lie algebra action $\varrho$ is called Hamiltonian if and
only if there is an element $J_0 \in C^1(\mathfrak g,\Cinf{M})$
such that $X_{\JN{\xi}} = X_\xi$ for all $\xi \in
\mathfrak g$, i.e. $i_{X_\xi} \omega= \d \JN{\xi}$. In this case
$\varrho(\xi)\cdot = \{\JN{\xi},{}\cdot{}\}$ and $J_0$ is said to
be a Hamiltonian for the action $\varrho$ (For applications in
physics where typically $\mathfrak g$ is the real Lie algebra
corresponding to a Lie group that acts on $M$ by symplectomorphisms
and where the generating vector fields $X_\xi$ are real-valued the
Hamiltonian $J_0$ is assumed to be real-valued, too.). In case
$J_0$ is equivariant with respect to the coadjoint representation
of $\mathfrak g$, i.e. $\{\JN{\xi},\JN{\eta}\}=
\JN{[\xi,\eta]}$ for all $\xi,\eta\in
\mathfrak g$ one calls $J_0$ a classical momentum mapping.

\begin{DEFINITION}
Let $\star$ be a $\mathfrak g$-invariant star product, then $J
= J_0 + J_+\in C^1(\mathfrak g,\Cinf{M})[[\nu]]$ with $J_0\in
C^1(\mathfrak g,\Cinf{M})$ and $J_+\in \nu C^1(\mathfrak
g,\Cinf{M})[[\nu]]$ is called a quantum Hamiltonian for the action
$\varrho$ in case
\begin{equation}\label{preqmmEq}
\varrho (\xi) = \frac{1}{\nu} \ad_\star (\Jbold{\xi})\quad
\textrm{for all}\quad\xi \in \mathfrak g.
\end{equation}
$J$ is called a quantum momentum mapping if in addition
\begin{equation}\label{qmmEq}
\frac{1}{\nu}\left(\Jbold{\xi}\star \Jbold{\eta}-
\Jbold{\eta}\star \Jbold{\xi}\right) = \Jbold{[\xi,\eta]}
\end{equation}
for all $\xi,\eta \in \mathfrak g$.
\end{DEFINITION}

Observe that the zeroth order in $\nu$ of (\ref{preqmmEq}) is
equivalent to $J_0$ being a Hamiltonian for $\varrho$ and that the
zeroth order in $\nu$ of (\ref{qmmEq}) just means equivariance of
this classical Hamiltonian with respect to the coadjoint action of
$\mathfrak g$ or equivalently that $J_0$ is a classical momentum
mapping. For Fedosov star products the fact that $J_0$ has to be a
classical Hamiltonian for $\varrho$ can also be seen directly from
Proposition \ref{quainnDerProp} as we have the following:

\begin{DEDUCTION}\label{qHamDed}
A $\mathfrak g$-invariant Fedosov star product for $(M,\omega)$
obtained from $(\nabla,\Omega,s)$ admits a quantum Hamiltonian if
and only if there is an element $J \in C^1(\mathfrak
g,\Cinf{M})[[\nu]]$ such that
\begin{equation}\label{preqmmFedEq}
\d \Jbold{\xi} = i_{X_\xi}(\omega+ \Omega)\quad \forall \xi \in
\mathfrak g \iff [i_{X_\xi}(\omega+\Omega)]=[0]\quad \forall \xi
\in \mathfrak g
\end{equation}
and from this equation $J$ is determined (in case it exists) up to
elements in $C^1(\mathfrak g,\mathbb C)[[\nu]]$.
\end{DEDUCTION}

\begin{REMARK}
Observe that the condition $H^1_{\mbox{\rm\tiny dR}}(M)=0$ is
obviously sufficient for the existence of a quantum Hamiltonian for
an arbitrary $\mathfrak g$-invariant star product $\star$ since
then any $\mathbb C[[\nu]]$-linear derivation of $\star$ is
quasi-inner. But for $\mathfrak g$-invariant Fedosov star products
$*$ the condition for the existence of a quantum Hamiltonian is
much weaker and more precise since only the cohomology classes of
very special closed one-forms have to vanish and not the complete
cohomology.
\end{REMARK}

Now recall the definition of a strongly invariant star product from
\cite{ArnCor83}:

\begin{DEFINITION}
Let $J_0$ be a classical momentum mapping for the action $\varrho$.
Then a $\mathfrak g$-invariant star product is called strongly
invariant if and only if $J=J_0$ defines a quantum Hamiltonian for
this action.
\end{DEFINITION}

Observe that the notion of strong invariance does not depend on the
chosen classical momentum mapping since every classical momentum
mapping is of the form $J_0 + b$ with $b\in Z^1_0(\mathfrak
g,\mathbb C)$ and hence every classical momentum mapping defines a
quantum Hamiltonian for $\varrho$ in case $J_0$ does. Moreover, in
the case of a strongly invariant star product $\star$ every
classical momentum mapping $J_0$ obviously yields a quantum
momentum mapping $J=J_0$ since $\frac{1}{\nu}\ad_\star
(\JN{\xi})\JN{\eta}= \{\JN{\xi},\JN{\eta}\} =\JN{[\xi,\eta]}$ for
all $\xi,\eta\in \mathfrak g$. As an immediate corollary of
Deduction \ref{qHamDed} we have:

\begin{COROLLARY}
Let $J_0$ be a classical momentum mapping for the action $\varrho$.
Then a $\mathfrak g$-invariant Fedosov star product $*$ obtained
from $(\nabla,\Omega,s)$ is strongly invariant if and only if
\begin{equation}\label{iXOmNul}
i_{X_\xi}\Omega = 0 \quad\textrm{for all}\quad\xi \in \mathfrak g.
\end{equation}
In this case every classical momentum mapping defines a quantum
momentum mapping for $*$.
\end{COROLLARY}
\begin{PROOF}
According to Deduction \ref{qHamDed} a classical momentum mapping
$J_0$ defines a quantum Hamiltonian for $*$ if and only if $\d
\JN{\xi} = i_{X_\xi}(\omega+\Omega)$ for all $\xi\in \mathfrak g$
but because of $\d \JN{\xi} = i_{X_\xi}\omega$ this is equivalent
to equation (\ref{iXOmNul}).
\end{PROOF}

Returning to the general case our next aim is to give a further
condition involving $\omega$, $\Omega$ and $X_\cdot$ which in
addition guarantees that a quantum Hamiltonian $J$ is in fact a
quantum momentum mapping.

\begin{PROPOSITION}\label{QmmProp}
Let $J$ be a quantum Hamiltonian for the Fedosov star product $*$
then $\lambda \in C^2(\mathfrak g,\Cinf{M})[[\nu]]$ defined by
\begin{equation}\label{lamDefEq}
\lambda(\xi,\eta):= \frac{1}{\nu}\left(\Jbold{\xi}* \Jbold{\eta}-
\Jbold{\eta} * \Jbold{\xi}\right) - \Jbold{[\xi,\eta]}
\end{equation}
lies in $C^2(\mathfrak g,\mathbb C)[[\nu]]$ and is an element of
$Z_0^2(\mathfrak g,\mathbb C)[[\nu]]$ which is explicitly given by
\begin{equation}\label{lamexplEq}
\lambda(\xi,\eta)=(\omega + \Omega)( X_\xi, X_\eta) -
\Jbold{[\xi,\eta]}
\end{equation}
and the cohomology class $[\lambda]\in H_0^2(\mathfrak g,\mathbb
C)[[\nu]]$ does not depend on the choice of $J$. Moreover quantum
momentum mappings exist if and only if $[\lambda]=[0]\in
H_0^2(\mathfrak g,\mathbb C)[[\nu]]$ and for every $a\in
C^1(\mathfrak g,\mathbb C)[[\nu]]$ such that $\delta_0 a = \lambda$
the element $J^a := J - a \in C^1(\mathfrak g,\Cinf{M})[[\nu]]$ is
a quantum momentum mapping for $*$. Finally, the quantum momentum
mapping (if it exists) is unique up to elements in $Z_0^1(\mathfrak
g,\mathbb C)[[\nu]]$, and hence we have uniqueness if and only if
$H_0^1(\mathfrak g,\mathbb C)=0$.
\end{PROPOSITION}
\begin{PROOF}
In fact all the statements of the proposition except for the
explicit shape of $\lambda$ hold for any $\mathfrak g$-invariant
star product $\star$ according to \cite[Prop. 6.3]{Xu98} and are
straightforward to prove. It thus remains to prove
(\ref{lamexplEq}) but this follows from the following computation
using equation (\ref{preqmmFedEq}):
\begin{eqnarray*}
\lefteqn{\lambda(\xi,\eta) + \Jbold{[\xi,\eta]}}\\
&=&
\frac{1}{\nu}\ad_*(\Jbold{\xi})\Jbold{\eta} = - \Lie_{X_\xi}
\Jbold{\eta} = - i_{X_\xi} \d \Jbold{\eta} = -
i_{X_\xi} i_{X_\eta}(\omega+\Omega) = (\omega +\Omega)( X_\xi,
X_\eta).
\end{eqnarray*}
\end{PROOF}

Again, for Fedosov star products the second condition for the
existence of a quantum momentum mapping can be formulated more
precisely than in the general case since the cocycle $\lambda$
whose cohomology class has to vanish to get a quantum momentum
mapping can be expressed explicitly in terms of $\omega$, $\Omega$
and $X_{\cdot}$. Obviously, supposing the existence of a classical
Hamiltonian for $\varrho$ the zeroth order of this condition is
equivalent to the existence of a classical momentum mapping.

Let us consider the important example of a semi-simple Lie algebra
$\mathfrak g$ in more detail:

\begin{EXAMPLE}
In case $\mathfrak g$ is semi-simple we have the following
properties: $[\mathfrak g, \mathfrak g]=\mathfrak g (\Rightarrow
H_0^1(\mathfrak g,\mathbb C)=0)$ and $H_0^2(\mathfrak g,\mathbb
C)=0$. But then $[\mathfrak g, \mathfrak g]=\mathfrak g$ implies
writing $\xi =\sum_{k \in I}[\zeta^{(k)},\eta^{(k)}]$ (the sum
ranges over a finite index set $I$) with $\zeta^{(k)},
\eta^{(k)}\in\mathfrak g$ and using the invariance of
$\omega+\Omega$ with respect to $X_{\zeta^{(k)}}$ and
$X_{\eta^{(k)}}$ that
\begin{equation*}
\begin{split}
i_{X_\xi}(\omega +\Omega) &= - \sum_{k
\in I}i_{[ X_{\zeta^{(k)}},X_{\eta^{(k)}}]}
(\omega +\Omega)\\&=-\sum_{k\in I}\Lie_{X_{\zeta^{(k)}}}
i_{X_{\eta^{(k)}}} (\omega +\Omega) =
\d\left(\sum_{k\in I}(\omega +\Omega)
(X_{\zeta^{(k)}},X_{\eta^{(k)}})\right)
\end{split}
\end{equation*}
and hence for all $\xi \in \mathfrak g$ there is a $\Jbold{\xi}\in
\Cinf{M}[[\nu]]$ such that $\d\Jbold{\xi} = i_{X_\xi}
(\omega +\Omega)$. Moreover, one can achieve that $J\in
C^1(\mathfrak g,\Cinf{M})[[\nu]]$ implying that $J$ defines a
quantum Hamiltonian for $*$ (e.g. fix a basis $\{e_i\}_{1\leq i\leq
\dim(\mathfrak g)}$ of $\mathfrak g$, write $e_i = \sum_{k
\in I_i}[\zeta^{(k)}_i, \eta^{(k)}_i]$, define $\Jbold{e_i}:=
\sum_{k\in I_i}(\omega +\Omega)
(X_{\zeta^{(k)}_i}, X_{\eta^{(k)}_i})$ such that $\d\Jbold{e_i}=
i_{X_{e_i}}(\omega +\Omega)$ holds according to the above
computation and extend $J$ to $\mathfrak g$ by linearity yielding
$J\in C^1(\mathfrak g,\Cinf{M})[[\nu]]$ with $\d\Jbold{\xi}=
i_{X_\xi}(\omega+\Omega)\,\forall
\xi\in\mathfrak g$.). This observation together with the statements
of Proposition \ref{QmmProp} and $H_0^1(\mathfrak g,\mathbb
C)=H_0^2(\mathfrak g,\mathbb C)=0$ implies that in this case there
is a unique quantum momentum mapping for every $\mathfrak
g$-invariant Fedosov star product.
\end{EXAMPLE}

Returning to the case of an arbitrary Lie algebra $\mathfrak g$ we
also have the following:

\begin{COROLLARY}
Let $*$ be a $\mathfrak g$-invariant Fedosov star product and
assume that there is a classical momentum mapping $J_0$ for the
action $\varrho$, then a quantum momentum mapping $J$ exists if and
only if there is an element $J_+\in
\nu C^1(\mathfrak g,\Cinf{M})[[\nu]]$ such that
\begin{equation}\label{qmmExClassExEq}
i_{X_\xi} \Omega = \d \JP{\xi}\quad\textrm{ and }\quad
\Omega(X_\xi,X_\eta) = (\delta_\varrho J_+)(\xi,\eta) \quad
\forall \xi,\eta \in \mathfrak g,
\end{equation}
and these equations determine $J_+$ up to elements of $\nu
Z^1_0(\mathfrak g,\mathbb C)[[\nu]]$.
\end{COROLLARY}
\begin{PROOF}
Assuming the existence of a classical momentum mapping it is
obvious that (\ref{preqmmFedEq}) and the equation
$\lambda(\xi,\eta)=0$ for all $\xi,\eta \in \mathfrak g$ reduce to
$i_{X_\xi}\Omega=\d\JP{\xi}$ and $\JP{[\xi,\eta]} =
\Omega(X_\xi,X_\eta)$ and it is straightforward to see that these
two equations are equivalent to (\ref{qmmExClassExEq}). The
statement about the ambiguity of $J_+$ is obvious from Proposition
\ref{QmmProp}.
\end{PROOF}

Observe that the condition for the existence of a quantum momentum
mapping for $\mathfrak g$-invariant Fedosov star products given in
the above corollary does not depend on the chosen classical
momentum mapping but only on $\Omega$ and $X_{\cdot}$. Moreover,
our result shows that the answer to the question whether existence
of a classical momentum mapping implies the existence of a quantum
momentum mapping posed in \cite{Xu98} in general is no if one
allows for star products whose characteristic class is different
from $\frac{1}{\nu}[\omega]$ since the conditions above involve the
two-form $\Omega$ that determines this class (cf. \cite{Neu99}) and
that has to be different from zero in this case. One can even
construct very simple examples where $\Omega$ is even exact and
hence the characteristic class is equal to $\frac{1}{\nu}[\omega]$
but nevertheless there exists no quantum momentum mapping.

\section*{Outlook and open Problems}
\label{ProbSec}
Let us conclude with a few remarks on our results and some possible
generalizations and questions that could be studied in the future:
\begin{enumerate}
\item
It should be possible to adapt our investigations to the case of
star products of Wick type on Semi-K\"ahler manifolds by imposing
additional conditions on the compatibility of the Lie algebra
action with the complex structure and due to the results of
\cite{Neu02a} such investigations would give a complete answer for
all such star products. These investigations will be subject of a
future project.
\item
A second possibility for generalizations could be to weaken the
conditions imposed on a quantum momentum mapping and to drop the
condition that $\frac{1}{\nu}\ad_\star(\Jbold{\xi})$ should equal
the Lie derivative with respect to $-X_\xi$ but to stick to the
condition of quantum covariance (\ref{qmmEq}) (which is reasonable
since this notion behaves properly with respect to equivalence
transformations of star products, which is not the case for the
notion of quantum momentum mappings considered in this paper) and
to demand that $\frac{1}{\nu}\ad_\star(\Jbold{\xi})=
-\Lie_{X_\xi} + O(\nu)$ is merely a deformation of the classical
Lie algebra action $\varrho$. Actually our results that establish a
strong relation between the characteristic class of the Fedosov
star product and the question of existence of a quantum momentum
mapping suggest that such a relation should also exist in general.
Maybe the fact that any star product is equivalent to a Fedosov
star product together with the results of the present paper can be
used to obtain results in this direction.
\end{enumerate}
\begin{small}

\end{small}
\end{document}